\documentclass[a4paper, 11pt, leqno]{amsart}
\usepackage[mathscr]{eucal}
\usepackage{amsmath}
\usepackage{amssymb}
\usepackage{amsfonts}
\usepackage{amsthm}
\theoremstyle{plain}
\newtheorem{theorem}{Theorem}[section]
\newtheorem{proposition}{Proposition}[section]
\newtheorem{corollary}{Corollary}[section]
\newtheorem{remark}{Remark}

\newtheorem{example}{Example}[section]

\newtheorem{problem}{Problem}[section]

\numberwithin{equation}{section}

\title{A Weierstrass type representation for minimal surfaces in Sol}
\author{Jun-ichi Inoguchi}
\thanks{The first named author is partially supported by Kakenhi 18540068}

\address{Department of Mathematics Education,
Utsunomiya University, \newline
Utsunomiya, 321-8505,
Japan}
\email{inoguchi@cc.utsunomiya-u.ac.jp}
\author{Sungwook Lee}
\address{Department of Mathematics, 
University of Southern Mississippi, 
Southern Hall, Box 5045, 
Hattiesburg, MS39406-5045 
U.S.A.}
\email{sunglee@usm.edu}

\dedicatory{Dedicated to professor Takeshi Sasaki on his 60th birthday}
\date{}
\begin{document}

\begin{abstract}
The normal Gauss map of a minimal surface in 
the model space $\mathrm{Sol}$ of solvegeometry
is a harmonic map with respect to a certain
singular 
Riemannian metric on the extended complex
plane. 
\end{abstract}

\keywords{Solvable Lie groups, minimal surfaces}

\subjclass[2000]{53A10, 53C15, 53C30}

\maketitle

\section{Introduction}

Since the discovery of holomorphic quadratic differential 
(called \textit{generalized Hopf differential} or
\textit{Abresch-Rosenberg differential}) for
CMC surfaces (constant mean curvature surfaces) in 3-dimensional homogeneous 
Riemannian manifolds with 
$4$-dimensional isometry group, 
global geometry of constant mean curvature surfaces 
in such spaces has been extensively studied \cite{AR}--\cite{AR2}. 

D.~A.~Berdinski{\u\i}
and
I.~A.~Ta{\u\i}manov \cite{BT} gave a representation formula for minimal
surfaces in 3-dimensional Lie groups in terms of spinors and Dirac operators.

The simply connected homogeneous Riemannian $3$-manifolds with 
\newline
\noindent
$4$-dimensional isometry group have structure of principal fiber bundle with
$1$-dimensional fiber and constant curvature base.
More explicitly, such homogeneous spaces are one of the following spaces;
the Heisenberg group $\mathrm{Nil}_3$,
the universal covering $\widetilde{\mathrm{SL}}_{2}\mathbb{R}$
of the special linear group equipped with naturally reductive metric,
the special unitary group $\mathrm{SU}(2)$ equipped with the Berger sphere
metric,
and reducible Riemannian symmetric space $S^{2}\times \mathbb{R}$,
$H^2\times \mathbb{R}$.

On the other hand, the model spaces of
Thurston's 3-dimensional
model geometries \cite{Thurston} are space forms, 
$\mathrm{Nil}$,
$\widetilde{\mathrm{SL}}_{2}\mathbb{R}$
with naturally reductive metric,
$S^{2}\times \mathbb{R}$,
$H^2\times \mathbb{R}$ and the space $\mathrm{Sol}$, 
the model space of solvegeometry.

Abresch and Rosenberg showed that 
the existence of generalized Hopf differential
in a simply connected Riemannian $3$-manifold
is equivalent to 
the property that the ambient space has at least $4$-dimensional
isometry group
\cite[Theorem 5]{AR2}.
Note that if the dimension of the isometry group
of a Riemannian $3$-manifold is greater than $3$, then 
the action of isometry group is transitive.

Thus for the space $\mathrm{Sol}$, one can not expect 
Abresch-Rosenberg type quadratic differential for 
CMC surfaces. 
Berdinski{\u\i}
and
Ta{\u\i}manov pointed out there are some difficulty to develop 
minimal surface geometry in $\mathrm{Sol}$ by using their 
representation formula and Dirac operators (see \cite[Remark 4]{BT}).

Thus, another approach for CMC surface geometry in $\mathrm{Sol}$ is
expected.

\vspace{0.2cm}

The space $\mathrm{Sol}$ belongs to the following two
parameter family of simply connected homogeneous Riemannian 3-manifolds;
$$
G(\mu_1,\mu_2)=(\mathbb{R}^3(x^1,x^2,x^3),g_{(\mu_1,\mu_2)}),
$$ 
with group structure
$$
(x^1,x^2,x^3)\cdot
(\tilde{x}^1,\tilde{x}^2,\tilde{x}^3)=
(x^1+e^{\mu_1x^3}\tilde{x}^1,
x^2+e^{\mu_2x^3}\tilde{x}^2,
x^3+\tilde{x}^3)
$$
and left invariant metric
$$
g_{(\mu_1,\mu_2)}=
e^{-2\mu_1x^3}(dx^{1})^{2}
+e^{-2\mu_2x^3}(dx^{2})^{2}+(dx^{3})^{2}.
$$
This family includes $\mathrm{Sol}=G(1,-1)$ as well as Euclidean 3-space 
$\mathbb{E}^3=G(0,0)$, hyperbolic 3-space $H^{3}=G(1,1)$ and
$H^{2}\times \mathbb{R}=G(0,1)$.

In this paper, we study 
the (normal)
Gauss map of minimal surfaces in $G(\mu_1,\mu_2)$. 
In particular, we shall show that 
the normal Gauss map of non-vertical minimal surfaces
is a harmonic map with respect to appropriate metric
if and only if $\mu_1^2=\mu_2^2$.

As a consequence, we shall give a Weierstrass-type representation formula 
for minimal surfaces in $\mathrm{Sol}$.
 
The results of this article were
partially reported at London Mathematical Society
Durham Conference ``Methods of Integrable Systems in Geometry" (August, 2006).

\section{Solvable Lie group}

In this paper, we study the following 
two-parameter family of homogeneous
Riemannian $3$-manifolds;
\begin{equation}\label{1.1}
\left\{(\mathbb{R}^{3}(x^{1},x^{2},x^{3}), g_{(\mu_1,\mu_2)})\
\vert \
(\mu_1,\mu_2)\in \mathbb{R}^{2}
\right \},
\end{equation}
where the metrics $g=g_{(\mu_1,\mu_2)}$ are
defined by
\begin{equation}
g_{(\mu_1,\mu_2)}:=e^{-2\mu_1x^3}(dx^{1})^{2}
+e^{-2\mu_2x^3}(dx^{2})^{2}+(dx^{3})^{2}.
\end{equation}
Each homogeneous space 
$(\mathbb{R}^3,g_{(\mu_1,\mu_2)})$
is realized as the following solvable matrix Lie group:
$$
G(\mu_1,\mu_2)=
\left\{\left(
\begin{array}{cccc}
1 & 0 & 0 & x^{3}\\
0 & e^{\mu_1x^3} & 0 & x^{1}\\
0 & 0 & e^{\mu_2x^3} & x^{2}\\
0 & 0 & 0 & 1
\end{array}
\right)
\
\Biggr
\vert
\
x^1,x^2,x^3
\in \mathbb{R}
\right
\}.
$$
The Lie algebra $\mathfrak{g}(\mu_1,\mu_2)$ is
given explicitly by
\begin{equation}
\mathfrak{g}(\mu_1,\mu_2)=
\left\{\left(
\begin{array}{cccc}
0 & 0 & 0 & y^{3}\\
0 & \mu_1y^3 & 0 & y^{1}\\
0 & 0 & \mu_2y^3 & y^{2}\\
0 & 0 & 0 & 1
\end{array}
\right)
\
\Biggr
\vert
\
y^1,y^2,y^3
\in \mathbb{R}
\right
\}.
\end{equation}
Then we can take the following
orthonormal basis $\{E_1,E_2,E_3\}$ 
of $\mathfrak{g}(\mu_1,\mu_2)$:
$$
E_{1}=
\left(
\begin{array}{cccc}
0 & 0 & 0 & 0\\
0 & 0 & 0 & 1\\
0 & 0 & 0 & 0\\
0 & 0 & 0 & 0
\end{array}
\right),
E_{2}=
\left(
\begin{array}{cccc}
0 & 0 & 0 & 0\\
0 & 0 & 0 & 0\\
0 & 0 & 0 & 1\\
0 & 0 & 0 & 0
\end{array}
\right),
E_{3}=
\left(
\begin{array}{cccc}
0 & 0 & 0 & 1\\
0 & \mu_1 & 0 & 0\\
0 & 0 & \mu_2 & 0\\
0 & 0 & 0 & 0
\end{array}
\right).
$$
Then the commutation relation of
$\mathfrak{g}$ is given by
$$
[E_1,E_2]=0,\
[E_2,E_3]=-\mu_{2}E_{2},\
[E_3,E_1]=\mu_{1}E_{1}.
$$
Left-translating the basis 
$\{E_1,E_2,E_3\}$, we obtain the following 
orthonormal frame field:
$$
e_{1}=e^{\mu_{1}x^3}\frac{\partial}{\partial x^{1}},
\
e_{2}=e^{\mu_{2}x^3}\frac{\partial}{\partial x^{2}},
\
e_{3}=\frac{\partial}{\partial x^{3}}.
$$
One can easily check that every $G(\mu_1,\mu_2)$ is 
a \textit{non-unimodular Lie group} except $\mu_1=\mu_2=0$.

The Levi-Civita connection $\nabla$ of 
$G(\mu_1,\mu_2)$ is described by
\begin{equation}\label{Levi-Civita} 
\begin{array}{ccc}
\nabla_{e_1}e_{1}=\mu_{1}e_{3},
& \nabla_{e_1}e_{2}=0,
& \nabla_{e_1}e_{3}=-\mu_{1}e_{1},
\\
\nabla_{e_2}e_{1}=0,\
& \nabla_{e_2}e_{2}=\mu_{2}e_{3},\
& \nabla_{e_2}e_{3}=-\mu_{2}e_{2},
\\
\nabla_{e_3}e_{1}=0,
& \nabla_{e_3}e_{2}=0,
& \nabla_{e_3}e_{3}=0.
\end{array}
\end{equation}
\begin{example}[Euclidean $3$-space]
{\rm
The Lie group $G(0,0)$ is
isomorphic and isometric to the Euclidean $3$-space
$\mathbb{E}^3=(\mathbb{R}^3,+)$. 
}
\end{example}

\begin{example}[Hyperbolic $3$-space]
{\rm
Take $\mu_1=\mu_2=c \not=0$. Then 
$G(c,c)$ is a warped product model of the
hyperbolic $3$-space:
$$
H^{3}(-c^2)=(\mathbb{R}^3(x^1,x^2,x^3),
e^{-2cx^3}\{(dx^{1})^{2}+(dx^{2})^{2}\}+(dx^{3})^{2}).
$$ 
}
\end{example}

\begin{example}[Riemannian product $H^2(-c^2)\times \mathbb{E}^1$]
{\rm
Take $(\mu_1,\mu_2)=(0,c)$ with $c\not=0$. Then
the resulting homogeneous space is
$\mathbb{R}^3$ with metric:
$$
(dx^{1})^{2}+e^{-2cx^3}(dx^{2})^{2}+(dx^{3})^{2}.
$$
Hence $G(0,c)$ is identified with
the Riemannian direct product of the Euclidean line
$\mathbb{E}^{1}(x^1)$ and the warped 
product model
$$
(\mathbb{R}^2(x^2,x^3),
e^{-2cx^{3}}(dx^{2})^{2}+(dx^{3})^{2}\
)
$$
of $H^{2}(-c^2)$.
Thus $G(0,c)$ is identified with 
$\mathbb{E}^1 \times
H^{2}(-c^2)$.
}
\end{example}  

\begin{example}[Solvmanifold]{\rm
The model space $\mathrm{Sol}$ of the $3$-dimensional
\textit{solvegeometry} \cite{Thurston} is
$G(1,-1)$. The Lie group $G(1,-1)$ is isomorphic to
the Minkowski motion group
$$
E(1,1):=\left\{
\left(
\begin{array}{ccc}
e^{x^3} & 0 & x^1 \\
0 & e^{-x^3} & x^2\\
0 & 0 & 1
\end{array}
\right)
\
\Biggr \vert
\
x^1,x^2,x^3
\in \mathbb{R}
\right\}.
$$
The full isometry group is $G(1,-1)$ itself.
}
\end{example}

\begin{example}\label{12plane}
{\rm
Since $[e_1,e_2]=0$, the distribution $D$ spanned by $e_1$
 and $e_2$ is involutive.
The maximal integral surface $M$ of $D$ through 
a point $(x^1_0,x^2_0,x^3_0)$ 
is the plane $x^3=x^3_0$.
One can see that $M$ is flat of constant mean
curvature $(\mu_1+\mu_2)/2$ (see (\ref{Levi-Civita}) 
). 

\begin{enumerate}
\item If $(\mu_1,\mu_2)=(0,0)$ then
$M$ is a totally geodesic plane.  
\item If $\mu_1=\mu_2=c\not=0$. Then $M$ is a horosphere in 
the hyperbolic $3$-space $H^3(-c^2)$.
\item If $\mu_1=-\mu_2\not=0$. Then $M$ is 
a non-totally geodesic minimal surface.
\end{enumerate}
} 
\end{example}

\section{Integral representation formula}

Let $M$ be a Riemann surface and
$(\mathfrak{D},z)$ be a simply connected 
coordinate region. The exterior derivative $d$ is decomposed as
$$
d=\partial+\bar{\partial},\
\partial=\frac{\partial}{\partial z}dz,\
\bar{\partial}=\frac{\partial}{\partial {\bar z}}d{\bar z},
$$
with respect to the conformal structure of $M$.
Take a triplet
$\{\omega^{1},\omega^{2},\omega^{3}\}$ of (1,0)-forms
which satisfies the following differential system:
\begin{eqnarray}
\bar{\partial}\omega^{i}&=& \mu_{i}\overline{\omega^i}\wedge \omega^{3},\ i=1,2;\label{HME1}
\\
\bar{\partial}\omega^{3}&=& \mu_{1}\omega^{1}\wedge \overline{\omega^1}+
\mu_{2}\omega^{2}\wedge \overline{\omega^2}.
\label{HME2}
\end{eqnarray}
\begin{proposition}[\cite{InoguchiSol}]\label{formula1}
Let $\{\omega^1,\omega^2,\omega^3\}$ be a solution to
(\ref{HME1})-(\ref{HME2}) on a simply connected coordinate region
$\mathfrak{D}$. Then
$$
\varphi(z,\bar{z})=2\int^{z}_{z_0}
\mathrm{Re}\>
\left(
e^{\mu_{1}x^{3}(z,\bar{z})}\cdot \omega^{1},
e^{\mu_{2}x^{3}(z,\bar{z})}\cdot \omega^{2},
\omega^{3}
\right)
$$
is a harmonic map of $\mathfrak{D}$
into $G(\mu_1,\mu_2)$.
\newline 
Conversely, any harmonic map of $\mathfrak{D}$ into $G(\mu_1,\mu_2)$ 
can be represented in this form. 
\end{proposition}
Equivalently, 
the resulting harmonic map
$\varphi(z,\bar{z})$ is defined by the
following data:
\begin{equation}
\omega^{1}=e^{-\mu_1x^3}x^1_{z}dz,\
\omega^{2}=e^{-\mu_1x^3}x^2_{z}dz,\
\omega^{3}=x^3_{z}dz,
\end{equation}
where the coefficient functions are
solutions to 
\begin{eqnarray}\label{HMEomega}
& x^i_{z\bar{z}}-\mu_{i}(x^{3}_{z}x^{i}_{\bar z}+x^{3}_{\bar z}x^{i}_{z})=0
\label{HMEomega12},
\ (i=1,2)\\
& x^3_{z\bar{z}}+\mu_{1}e^{-2\mu_{1}x^3}x^{1}_{z}x^{1}_{\bar z}
+\mu_{2}e^{-2\mu_{2}x^3}x^{2}_{z}x^{2}_{\bar z}=0
\label{HMEomega3}.
\end{eqnarray}

\begin{corollary}[\cite{InoguchiSol}]
Let $\{\omega^1,\omega^2,\omega^3\}$ be a solution to
\begin{eqnarray}
& \bar{\partial}\omega^{i}= \mu_{i}\overline{\omega^i}\wedge \omega^{3},\ i=1,2;\\
& \omega^{1}\otimes \omega^{1}+
\omega^{2}\otimes \omega^{2}+
\omega^{3}\otimes \omega^{3}=0
\end{eqnarray}
on a simply connected coordinate region
$\mathfrak{D}$. Then
$$
\varphi(z,\bar{z})=2\int^{z}_{z_0}
\mathrm{Re}\>
\left(
e^{\mu_{1}x^{3}(z,\bar{z})}\cdot \omega^{1},
e^{\mu_{2}x^{3}(z,\bar{z})}\cdot \omega^{2},
\omega^{3}
\right)
$$
is a weakly conformal harmonic map
of $\mathfrak{D}$ into $G(\mu_1,\mu_2)$.
Moreover $\varphi(z,\bar{z})$ is a minimal immersion if
and only if 
$$
\omega^{1}\otimes \overline{\omega^{1}}+
\omega^{2}\otimes \overline{\omega^{2}}+
\omega^{3}\otimes \overline{\omega^{3}}\not=0.
$$
\end{corollary}

\section{The normal Gauss map}

Let $\varphi:M\to G(\mu_1,\mu_2)$ be a conformal 
immersion.
Take a unit normal vector 
field $N$ along $\varphi$.
Then, by the left translation we obtain the following smooth map:
$$
\psi:=dL_{\varphi}^{-1}\cdot N:M\to S^{2}\subset \mathfrak{g}(\mu_1,\mu_2).
$$
The resulting map $\psi$ takes value in the unit $2$-sphere
$S^2$ in the Lie algebra $\mathfrak{g}(\mu_1,\mu_2)$.
Here, 
via the orthonormal basis $\{E_1,E_2,E_3\}$,
we identify $\mathfrak{g}(\mu_1,\mu_2)$ with Euclidean $3$-space
$\mathbb{E}^3(u^1,u^2,u^3)$.

The smooth map $\psi$ is called the \textit{normal Gauss map}
of $\varphi$.

Let $\varphi:\mathfrak{D}\to G(\mu_1,\mu_2)$ be a
weakly conformal
harmonic map of a simply connected Riemann surface $\mathfrak{D}$
determined by the data $(\omega^1,\omega^2,\omega^3)$.
Express the data as $\omega^i=\phi^{i}dz$.
Then the induced metric $I$ of $\varphi$ is 
$$
I=2(\sum_{i=1}^{3}|\phi^{i}|^{2})dzd{\bar z}. 
$$
Moreover these three coefficient functions satisfy
$$
\frac{\partial 
\phi^3}{\partial{\bar z}}=-\sum_{i=1}^{2}\mu_{i}|\phi^i|^{2},\ \
\frac{\partial
\phi^{i}}{\partial{\bar z}}=\mu_{i}\>
\overline{\phi^i}\>
\phi^{3}, \ i=1,2,
$$

\begin{equation}\label{nullcurve}
(\phi^1)^2+(\phi^2)^2+(\phi^3)^{2}=0.
\end{equation}
The harmonic map $\varphi$ is 
a minimal immersion if and 
only if 
\begin{equation}\label{nonsingular}
|\phi^1|^2+|\phi^2|^2+|\phi^3|^{2}\not=0.
\end{equation}

Here we would like to remark that
$\phi^3$ is identically 
zero if and only if $\varphi$ is a 
\textit{vertical plane}
$x^3=\mathrm{constant}$. 
(See example \ref{12plane}).
As we saw in example \ref{12plane},
the vertical plane 
$\varphi$ is minimal if and only if 
$\mu_1+\mu_2=0$.

Hereafter we assume that $\phi^3$ is not identically zero.
Then
we can introduce two mappings $f$ and $g$ by
\begin{equation}\label{data-fg}
f:=\phi^1-\sqrt{-1}\phi^2,\
\
g:=\frac{\phi^3}{\phi^1-\sqrt{-1}\phi^2}.
\end{equation}
By definition, $f$ and $g$ take values in the extended
complex plane $\overline{\mathbb{C}}=\mathbb{C}\cup
\{\infty\}$.
Using these two 
$\overline{\mathbb{C}}$-valued functions, $\varphi$ is rewritten as
$$
\varphi(z,\bar{z})=2\int^{z}_{z_0}\mathrm{Re}
\left(
e^{\mu_1x^3}\frac{1}{2}f(1-g^2),
e^{\mu_2x^3}\frac{\sqrt{-1}}{2}f(1+g^2),fg
\right)dz.
$$
The normal Gauss map
is computed as
$$
\psi(z,\bar{z})=\frac{1}{1+|g|^{2}}\left(
2\mathrm{Re}\>(g) E_{1}+
2\mathrm{Im}\>(g)E_{2}+
(|g|^2-1)E_{3} 
\right).
$$

Under the stereographic projection 
$\mathcal{P}:S^{2}\setminus \{\infty\}\subset \mathfrak{g}(\mu_1,\mu_2)\to
\mathbb{C}:=\mathbb{R}E_{1}+\mathbb{R}E_2$,
the map $\psi$
is identified with the $\overline{\mathbb{C}}$-valued
function $g$.
Based on this fundamental observation, we call the function $g$
the \textit{normal Gauss map} of $\varphi$.
The harmonicity together with the integrability
(\ref{HMEomega12})--(\ref{HMEomega3}) are equivalent to
the following system for $f$ and $g$:

\begin{eqnarray}
 \frac{\partial f}{\partial \bar{z}}&=&\frac{1}{2}|f|^{2}g\{
\mu_{1}(1-\bar{g}^2)
-\mu_{2}(1+\bar{g}^{2})
\}\label{f-eq}
,\\
\frac{\partial g}{\partial \bar{z}}&=&-\frac
{1}{4}
\{
\mu_{1}(1+g^2)(1-\bar{g}^2)
+\mu_2
(1-g^2)(1+\bar{g}^2)\}
\bar{f}.
\label{g-eq}
\end{eqnarray}

\begin{theorem}[\cite{InoguchiSol2}]\label{formula2}
Let $f$ and $g$ be a $\overline{\mathbb{C}}$-valued functions 
which are solutions to the system{\rm:}
(\ref{f-eq})--(\ref{g-eq}). Then 
\begin{equation}
\varphi(z,\bar{z})=2\int^{z}_{z_0}\mathrm{Re}
\left(
e^{\mu_1x^3}\frac{1}{2}f(1-g^2),
e^{\mu_2x^3}\frac{\sqrt{-1}}{2}f(1+g^2),fg
\right)dz
\end{equation}
is a weakly conformal harmonic map
of $\mathfrak{D}$ into $G(\mu_1,\mu_2)$.
\end{theorem}

\begin{example}
{\rm
Assume that $\mu_1\not=0$.
Take the following two $\overline{\mathbb{C}}$-valued 
functions: 
$$
f=\frac{\sqrt{-1}}{\mu_{1}(z+\bar{z})},
\
g=-\sqrt{-1}.
$$
Then $f$ and $g$ are solutions to
(\ref{f-eq})--(\ref{g-eq}). 
By the integral representation formula, we
can see that
the minimal surface determined by the data $(f,g)$ 
is a plane $x^2=\mathrm{constant}$.
Note that this plane is totally geodesic in $G(1,-1)$.
}
\end{example}

From (\ref{f-eq})--(\ref{g-eq}), we can eliminate 
$f$ and deduce the following PDE for $g$.

\begin{align}\label{eq:harm}
&{} g_{z\bar z}-\frac{2g\{\mu_1(1-{\bar g}^2)-\mu_2(1+{\bar
g}^2)\}g_zg_{\bar z}}{\mu_1(1+g^2)(1-{\bar
g}^2)+\mu_2(1-g^2)(1+{\bar g}^2)}
\\
\notag
&+\frac{4{\bar g}(1-g^4)(\mu_1^2-\mu_2^2)|g_{\bar z}|^2}{(\mu_1^2+\mu_2^2)|1-g^4|^2+\mu_1\mu_2\{(1+g^2)^2(1-{\bar
           g}^2)^2+(1+{\bar g}^2)^2(1-g^2)^2\}}\\
\notag
&=0.
\end{align}

\begin{theorem}
The equation \eqref{eq:harm} is the harmonic map equation for a map
$g: \mathfrak{D}\longrightarrow\overline{\mathbb C}(w,\bar w)$ if and only
if $\mu_1^2=\mu_2^2$.
\begin{enumerate}
\item
If $\mu_1=\mu_2\ne 0$, then the equation \eqref{eq:harm} becomes
\begin{equation}
\label{eq:harm2} 
\frac{\partial^2 g}{\partial z\partial\bar
z}+\frac{2|g|^2\bar g}{1-|g|^4}\frac{\partial g}{\partial
z}\frac{\partial g}{\partial\bar z}=0.
\end{equation}
The differential equation \eqref{eq:harm2} is the harmonic map equation for a map
$g$ from $\mathfrak{D}$ into $\left(\bar {\mathbb C}(w,\bar w),\frac{dwd\bar
w}{|1-|w|^4|}\right)$. The singular metric $\frac{dwd\bar
w}{|1-|w|^4|}$ is called the Kokubu metric {\rm(\cite{AA}, \cite{Kokubu})}. 
\item
If $\mu_1=-\mu_2\ne
0$, then \eqref{eq:harm} becomes
\begin{equation}
\label{eq:harm3} \frac{\partial^2 g}{\partial z\partial\bar
z}-\frac{2g}{g^2-{\bar g}^2}\frac{\partial g}{\partial
z}\frac{\partial g}{\partial\bar z}=0.
\end{equation}
The differential equation
\eqref{eq:harm3} is the harmonic map equation for a map $g$ from $\mathfrak{D}$
into $\left(\overline{\mathbb C}(w,\bar w),\frac{dwd\bar w}{|w^2-{\bar
w}^2|}\right)$.
\end{enumerate}
\end{theorem}
\noindent
\textit{Proof.} 
Consider a possibly singular Riemannian metric $\lambda^{2} dwd\bar{w}$ on 
the extended complex plane $\overline{\mathbb{C}}(w,\bar{w})$.
Denote by $\varGamma^{w}_{ww}$ 
the Christoffel symbol of the metric 
with respect to $(w,\bar{w})$. 
Then for a map $g: M\longrightarrow\bar{\mathbb
C}(w,\bar w)$, the \textit{tension field}
$\tau(g)$ of $g$ is given by
\begin{equation}
\label{eq:tension}
\tau(g)=4\lambda^{-2}\left(
g_{z\bar z}+\varGamma^w_{ww}g_zg_{\bar z}
\right).
\end{equation}
By comparing the equations
\eqref{eq:harm} and $\tau(g)=0$, one can readily see that
\eqref{eq:harm} is harmonic map equation if and only if
$\mu_1^2=\mu_2^2$.

In order to find a suitable metric on $\bar{\mathbb C}(w,\bar w)$
with which \eqref{eq:harm} is harmonic map equation, one simply
needs to solve the first order PDE
$$
\left\{\begin{aligned} 
\varGamma^w_{ww}&=\frac{2|w|^2\bar w}{1-|w|^4} &{\rm if}\ \mu_1&=\mu_2\ne 0,\\
\varGamma^w_{ww}&=-\frac{2w}{w^2-{\bar w}^2} &{\rm if}\
\mu_1&=-\mu_2\ne 0, \end{aligned}\right.
$$
whose solutions are
$\lambda^{2}=1/|1-|w|^4|$
and
$\lambda^{2}=1/|w^2-{\bar w}^2|$,
respectively. $\blacksquare$

\begin{corollary}
Let $g:\mathfrak{D}\to \left(\overline{\mathbb C}(w,\bar w),\frac{dwd\bar w}{|w^2-{\bar
w}^2|}\right)$ be a harmonic map.
Define a function $f$ on $\mathfrak{D}$ by
$$
f=\frac{2\bar{g}_z}{g^2-\bar{g}^2}.
$$
Then
$$
\varphi(z,\bar{z})=2\int^{z}_{z_0}\mathrm{Re}
\left(
e^{x^3}\frac{1}{2}f(1-g^2),
e^{-x^3}\frac{\sqrt{-1}}{2}f(1+g^2),fg
\right)dz
$$
is a weakly conformal harmonic map
of $\mathfrak{D}$ into $\mathrm{Sol}$.
\end{corollary}
\begin{remark}{\rm
Direct computation shows the
following formulas:
\begin{enumerate}
\item The sectional curvature of $(\overline{\mathbb{C}}(w,\overline{w}),dwd\bar{w}/|1-|w|^{4}|)$ is
\newline $-8|w|^{2}/|1-|w|^{4}|$.
\item
 The sectional curvature of $(\overline{\mathbb{C}}(w,\overline{w})
,dwd\bar{w}/|w^2-{\bar w}^{2}|)$ is
\newline
$-8|w|^{2}/|w^{2}-\bar{w}^2|$.
\end{enumerate}
}
\end{remark}
\begin{remark}
{\rm The normal Gauss map of a non-vertical minimal
surface in the Heisenberg group is a harmonic map into
the hyperbolic $2$-space. See \cite{InoguchiNil}.
}
\end{remark}

Aiyama and Akutagawa \cite{AA}
studied the Dirichlet problem at infinity for
proper harmonic maps from the unit disc to the extended complex plane
equipped with the Kokubu metric.
To close this paper we propose the following probelm:
\begin{problem}
Study Dirichlet problem at infinity for
harmonic maps into the extended complex plane with 
metric $dwd{\bar w}/|w^2-{\bar
w}^2|$ and apply it for the construction of
minimal surfaces in $\mathrm{Sol}$.
\end{problem}

\end{document}